\documentclass[a4paper,11pt]{article}
\usepackage{array}
\usepackage{theorem}
\usepackage{amsmath,amscd,amssymb} 
\usepackage{latexsym}
\usepackage{epsfig}
\theorembodyfont{\sl}

\newcommand{\CC}{\mathbb C}

\newcommand{\FF}{\mathbb F}

\newcommand{\HH}{\mathbb H}

\newcommand{\NN}{\mathbb N}

\newcommand{\QQ}{\mathbb Q}

\newcommand{\cA}{\mathcal A}

\newcommand{\cS}{\mathcal S}

\renewcommand{\Bar}{\overline}

\newcommand{\To}{\longrightarrow}

\newcommand{\GL}{\mathop{\mathrm {GL}}\nolimits}

\renewcommand{\Im}{\mathop{\mathrm {Im}}\nolimits}

\newcommand{\gothM}{\mathfrak M}

\newcommand{\cusp}{\mathop{\mathrm {cusp}}\nolimits}

\begin{document}
\title{A supersingular coincidence}
 \author{G.K. Sankaran}

\maketitle

The list of fifteen primes
$$
\cS=\{2,3,5,7,11,13,17,19,23,29,31,41,47,59,71\}
$$
known as the \emph{supersingular primes} ({\tt https://oeis.org/A002267}) appears
in several different contexts. Here are five of them.
\begin{enumerate}
\item $p\in \cS$ if and only if $p$ divides the order of the Monster sporadic
simple group.
\item $p\in \cS$ if and only if $g(X_0(p)^+)=0$, where $X_0(p)^+=\HH/\Gamma_0(p)^+$ is the
modular curve associated with the group $\Gamma_0(p)^+<\GL(2,\QQ)$
\item $p\in \cS$ if and only if the supersingular values of the $j$-invariant all lie in
  $\FF_p$. 
\item $p\in\cS$ if and only if the space $J^{\cusp}_{2,p}$ of Jacobi
  cusp forms of weight~$2$ and index~$p$
is of dimension~$0$.
\item If $p\not\in\cS$ then the moduli space $\cA_p$ of complex abelian surfaces with a
polarisation of type $(1,p)$ is of general type.
\end{enumerate}

Of these, (1)--(3) are described in~\cite{ogg}, where the equivalence
of the conditions in (2) and (3) is proved. In~\cite{ogg} a prize (a
bottle of Jack Daniels) is offered for an explanation of why the
condition in (1) is equivalent to those in (2) and (3): it is still
unclaimed.

This note is primarily about (5). The proof that $\cA_p$ is of general type
for $p\not\in\cS$ is due to Erdenberger~\cite{erdenberger}, and specialists
in moduli of abelian surfaces are occasionally asked to explain the
apparent coincidence~\cite{hemckay}. In fact the answer consists of a series of
well-known facts, but because they are not all well known to the same
people, the question continues to recur. The purpose of this note is
to set the answer out clearly.

\section{Moduli of abelian surfaces}\label{sect:moduli}

An abelian surface equipped with a polarisation of type $(1,d)$ (for
$d\in\NN$) may be thought of as a complex torus $\CC^2/\Lambda$, where
$\Lambda\subset \CC^2$ is the subgroup (lattice) generated by the
columns of $\Omega=(I_d,\tau)$ for
\[
I_d:=\begin{pmatrix} 1 & 0 \\ 0 & d
\end{pmatrix}, \qquad \tau = \begin{pmatrix}\tau_1 & \tau_2\\ \tau_2& \tau_3
\end{pmatrix} \in \HH_2=\{\tau={}^t\tau\in M_{2\times
2}(\CC)\mid \Im \tau >0\}.
\]
The \emph{paramodular group}
\[
\Gamma_d=\left\{\gamma\in \GL(4,\QQ)\mid {}^t\gamma \begin{pmatrix}0
& I_d\\ -I_d & 0
\end{pmatrix}\gamma = \begin{pmatrix}0
& I_d\\ -I_d & 0
\end{pmatrix}\right\}
\]
acts on the Siegel upper half-plane $\HH_2$ by fractional linear
transformations
\[
\begin{pmatrix}A&B\\ C&D
\end{pmatrix}\colon \tau \To (A\tau+B)(C\tau+D)^{-1}.
\]
This group action is properly discontinuous and the quotient
$\cA_d:=\HH_2/\Gamma_d$ is a coarse moduli space for $(1,d)$-polarised
abelian varieties. It is a quasi-projective variety, and one may ask
for its Kodaira dimension, or more precisely, for the Kodaira
dimension $\kappa(Y_d)$ of a desingularisation $Y_d$ of a projective
compactification $\Bar\cA_d$.  For more on this and related
spaces, see~\cite{HKW}.

In practice one expects that $\cA_d$ is of general type,
i.e.\ $\kappa(Y_d)=3$, except for some small values of $d$. Very
loosely, this is because $k$-fold differential forms on $Y_d$
correspond to suitable modular forms of weight $3k$ for $\Gamma_d$,
and these become abundant as $d$ grows at least for $k$ sufficiently
divisible. However, not every modular form of weight $3k$ will do:
obstructions come from the boundary $\Bar\cA_d\setminus \cA_d$ and
from the branching of $\HH_2\to \cA_d$.

The obstructions at the boundary may be overcome by using the
\emph{low-weight cusp form trick}~\cite{egloffstein}: if we can find a
cusp form $f_2$ of weight~$2$ for $\Gamma_d$ then we may consider
modular forms $f$ of weight~$3k$ of the form $f=f_2^kf_k$, where $f_k$
is a modular form of weight~$k$. These are also abundant, if $d$ and
$k$ are large enough, and because they vanish to high order at the
boundary, the associated differential forms extend.

The branching behaviour has to be analysed separately, and it depends
on the factorisation of $d$. For that reason much work in this
direction has concentrated, for simplicity, on the case $d=p$ prime.
The case $d=p^2$ has some simplifying features and was treated
in~\cite{og} and~\cite{gs}.

By this method it was shown in~\cite{ap} that $\cA_p$ is of general type
for $p>173$. Because of the inefficient compactification used there,
the effective constraint on $p$ came from the branching, so all that
was necessary was to verify that a weight~$2$ cusp form exists for all
$p>173$. Such a form may be obtained by lifting a Jacobi cusp form of
weight~$2$ and index~$p$ according to Gritsenko
(\cite[Theorem~3]{egloffstein}, \cite{gIMRN}). The dimension of the
space of Jacobi cusp forms is computed in~\cite{ez,sz} and in this
case it takes the form \cite{egloffstein,ap}
\[
\dim J^{\cusp}_{2,p}=\sum_{j=1}^p\left\lfloor\frac{1+j}{6}\right\rfloor-\delta_6(j)-\left\lfloor \frac{j^2}{4p}\right\rfloor
\]
where $\delta_6(j)=1$ if $6|j$ and $0$ otherwise. This is positive for
all $p>173$.

Erdenberger~\cite{erdenberger} found a better compactification and was
able to reduce the condition imposed by the branching from $p>173$ to
$p\ge 37$, so that the existence of the Jacobi form becomes the
effective constraint. Then it is easy to compute from the formula
above that $p\in \cS$ exactly when no Jacobi cusp form of weight~$2$
and index~$p$ exists, i.e.\ when the condition in (4) holds.

It is not necessarily to be expected that $\cA_p$ is of general type
exactly when $p\not\in \cS$. The method of proof of \cite{erdenberger}
fails for $p\in\cS$, as we shall see, and it is known that $\cA_d$ is
unirational (so in particular not of general type) for some small
values of $d$, including all primes $p\le 11$: see~\cite{gp}.
However, if $p\ge 13$, nothing currently excludes the possibility that
$\cA_p$ is of general type.

On the other hand, $\cA_p$ being unirational, other than in the known
cases $p\le 11$, \emph{is} excluded. Gritsenko~\cite{gIMRN} showed
that $\cA_d$ has non-negative Kodaira dimension, so is not uniruled,
for all $d\ge 13$, prime or not, except possibly for $d=14,
15,16,18,20,24,30, 36$. Of these, the cases $d=14,16,18, 20$ have
since been settled in~\cite{gp} (such $\cA_d$ are in fact
unirational) and only for $d=15,24,30,36$ is nothing known about the
Kodaira dimension of $\cA_d$.

\section{Modular forms}\label{sect:modular}

Since we have now established a connection between (4) and (5), to
achieve a moderately satisfactory explanation of the apparently
coincidental appearance of $\cS$ in (5) we should show, without direct
computation, that the conditions in (2) and (4) are equivalent. (A
fully satisfactory explanation would also involve (1): this we
are not able to give.) This is well known among specialists in
Jacobi forms, and follows easily from a small part of~\cite{sz}.

It is shown in \cite{sz} that the space $J_{k,d}$ of Jacobi forms of
weight $k$ and index $d$ is isomorphic (even as a Hecke module) to a
certain subspace $\gothM^-_{2k-2}(d)$ of the space $M_{2k-2}(d)$ of
modular forms of weight $2k-2$ for $\Gamma_0(d)$. This subspace is
defined by $\gothM^-_{2k-2}(d)=M^-_{2k-2}(d)\cap\gothM_{2k-2}(d)$,
where $M^-_{2k-2}(d)$ is the space of weight~$2k-2$ modular forms for
$\Gamma_0(d)$ that satisfy an extra condition on the behaviour under
the Fricke involution $w\colon \tau \mapsto\frac{-1}{d\tau}$, namely
\[
f(\frac{-1}{d\tau})=(-1)^k d^{k-1}\tau^{2k-2}f(\tau).
\]
In our case ($k=2$ and $d=p$) this is equivalent to saying that
$f$ is a modular form of weight $2$ for the group
$\Gamma_0(p)^+<\GL(2,\QQ)$ generated by $\Gamma_0(p)$ and
$w=\begin{pmatrix}0&1\\ -d&0
\end{pmatrix}$. See, for example, the definition of automorphic form
in~\cite[Chapter~2]{Shi}.  So if there is a weight~$2$, index~$p$
Jacobi form, then $\Gamma_0(p)^+$ has a weight~$2$ modular form.

Conversely, inspecting the definition of $\gothM$ in~\cite{sz}
we find that there are no other conditions for $p$ prime: simply
$\gothM_{2k-2}(p)=M_{2k-2}(p)$. This can be seen at once, for
instance, from \cite[Equation~(4), p.~116]{sz}, since
$\gothM_2(1)\subset M_2(1)=0$. In other words, the space of Jacobi
forms in this case is isomorphic exactly to the space of weight~$2$
modular forms for $\Gamma_0(p)^+$. Moreover, the isomorphism respects
cusp forms: see \cite[Theorem~5]{sz}.
\par
We remark that for $k=2$ and $p$ square-free (in particular for $p$
prime) there are no Eisenstein series, so the condition (3) is
equivalent to the same statement but with $J^{\cusp}_{2,p}$ replaced
by $J_{2,p}$.
\par
However, Ogg~\cite{ogg} shows that the modular curve $X_0(p)^+$
corresponding to $\Gamma_0(p)^+$ is of genus~$0$ precisely for
$p\in\cS$, i.e.\ he shows~(2). One can compute the dimension of the
space of weight~$2$ forms from the formulae given
in \cite[Theorem~2.23]{Shi}: it is $g+m-1$, where $g$ is the genus
of $X_0(p)^+$ and $m$ is the number of cusps. Because $p$ is prime,
the curve $X_0(p)$ has two cusps, which are interchanged by the Fricke
involution; so $m=1$, and so the space of modular forms for
$\Gamma_0(p)^+$ has dimension $g$ (i.e.\ they are all cusp forms, as
one should also expect from the remark above). So
for $p\in\cS$, there can be no weight~$2$, index~$p$ Jacobi forms; so
we definitely cannot prove that $\cA_p$ is of general type result by the methods
of~\cite{ap} and~\cite{erdenberger} for any $p\in\cS$.

G.K. Sankaran,\\ Department of Mathematical Sciences,\\ University of
Bath,\\ Bath BA2 7AY,\\ England\\ \\ {\tt G.K.Sankaran@bath.ac.uk}

\end{document}